
\mag=1200
\hsize=6.5 true in
\vsize=8.7 true in
  \baselineskip=15pt
\vglue 1.5 cm
\def\fnote#1{\footnote{$^{#1}$}}
\count1=\number\year
\advance\count1 by -1900
\def\abc{\number\month/\number\day/\number\count1}
\footline={\hss\tenrm -\folio-\ \ \abc \hss}
        


  \def\qed{$\rlap{$\sqcap$}\sqcup$}
\def\prim{^\prime }



\font\tengothic=eufm10
\font\sevengothic=eufm7
\newfam\gothicfam
      \textfont\gothicfam=\tengothic
      \scriptfont\gothicfam=\sevengothic
\def\goth#1{{\fam\gothicfam #1}}


   \font\tenmsb=msbm10              \font\sevenmsb=msbm7
\newfam\msbfam
      \textfont\msbfam=\tenmsb
      \scriptfont\msbfam=\sevenmsb
\def\Bbb#1{{\fam\msbfam #1}}
\def\PP#1{{\Bbb P}^{#1}}


\def\move-in{\parshape=1.75true in 5true in}

\def\ref#1{[{\bf #1}] }

\hyphenation {Castel-nuovo}


\medskip
\vglue 12 pt
\centerline{{\bf ON THE HILBERT FUNCTION OF FAT POINTS }}
\centerline{{\bf ON A RATIONAL NORMAL CUBIC}\fnote{*}{This work has been
supported by MURST funds and by the CNR group "Rami analitici e
sistemi lineari".}}
\medskip
\medskip
\centerline {M.V.Catalisano , A.Gimigliano}
\medskip
\medskip
\medskip
\medskip
\medskip
\medskip
\medskip
\medskip
\noindent {\bf Abstract: } In this paper we find an algorithm which computes
the
Hilbert function of schemes $Z$ of "fat points" in $\PP3$
whose support lies on a rational
normal cubic curve $C$.  The algorithm shows that the maximality of the Hilbert
function
 in degree $t$ is related to the existence of fixed curves (either $C$ itself
or
one of its secant lines) for the linear system of surfaces of degree $t$
containing $Z$.
\medskip
\medskip
\medskip
\medskip
\medskip
\noindent {\bf Introduction}
\par
\medskip
The aim of this paper is to consider linear systems $J_t$
defined by particular schemes of fat
points, where with "fat points" we mean $0$-dimensional
 schemes defined by homogeneous
 ideals of type
\par
\noindent(*) $\qquad \qquad \qquad \qquad \qquad
J=\oplus _{t\geq 0} J_t ={\goth p}^{m_1} _1 \cap \ldots
\cap {\goth p}^{m_s} _s\  $
\par
\noindent where each  ${\goth p}_{i}$ is the homogeneous ideal in
${\bf R} = k[x_{0},\ldots ,x_{r}]$ of a point
$P_{i} \in  \PP {r} = \PP {r} _k$  ($k$
 being an algebraically closed field of characteristic 0), and the
$m_i$'s
 are non negative integers. We will denote a scheme of fat points by
$Z = (P_1,...,P_s;m_1,...,m_s)$.
\par
In [{\bf 3}], a bound for the regularity  of the linear systems of
type $J_{t}$ is given when the points $P_{i} \in  \PP {r}$ are in
(linear) generic position (i.e. no three on a line, no four on a plane,
etc.).  It
turns out that the "worst" case for $J_{t}$ is when the points $P_{i}$
lie
on a rational normal curve (see also \ref{7}).
\par
This leads to the following conjecture:
\par
\medskip
\noindent {\bf Conjecture:} {\it Let} $J \subseteq  {\bf R}$
{\it be an ideal
of fat points in} $\PP {r} ${\it (i.e. J is as in} (*){\it ),  and let}
$H({\bf R}/J,t)$ {\it denote the Hilbert function of} ${\bf R}/J$.
{\it Then, if the points} $P_{i}$ {\it are in linear generic position,}
$\forall t\in  {\Bbb N}$  {\it we
have that} $H({\bf R}/J,t) \ge  H({\bf R}/I,t)$ {\it , where} $I$ {\it is an
ideal of type} (*), {\it with the same multiplicities }$m_i$ {\it as } $J $
{\it and whose support is given by points on a rational normal curve}
$C_{r}\subseteq  \PP {r}${\it.  Moreover the  value} $H({\bf R}/I,t)$
{\it
does not depend on the choice of the} s {\it points on} $C_{r}${\it .}
\par
\medskip
In this paper we analyze the case $r = 3$, and we show (via Theorem 2.2)
that there is an algorithm which computes $H({\bf R}/I,t)$ for ideals $I$
as above (i.e. for fat points on a cubic curve $C$). The algorithm
will only depend on the data $s,m_1,...,m_s$, thus showing that the
Hilbert
function  does not depend on the position of the (distinct)
points on the curve.
\par
It will also turn out that $H({\bf R}/I,t)$ has its maximal value (i.e.
the fat
points impose independent conditions to surfaces of degree $t$), if and
only if for every $P_i$ with $m_i > 0$ the linear system
$(I:{\goth p}_i)_t$ has
neither $C$ nor any line $P_iP_j$ as fixed locus (Corollary 2.3).
\par
\medskip
The paper is divided as follows: the first section is devoted to studying
the following question: which numerical ("Bezout-like") conditions
imply that a multiple of a curve (a line $P_iP_j$ or the curve $C$ in our
case) is a fixed locus for the linear system
$I_t$ ? In this section we also consider whether the numerical conditions
that we find are necessary and we
compute the Hilbert function of all the multiples of $C$, i.e.
of the ideals $I_C ^n$.
\par
In $\S 2$ we state the main result and describe it, while $\S 3$ is
dedicated to several lemmata which will be used in $\S 4$ to prove
Theorem 2.2.
\par
 We would like to thank C.H.Walter for some useful talks we had during
 the Workshop organized by the "Europroj" group on points (Nice '93),
when the work on this paper started.
\par
\medskip
\medskip
\medskip

\noindent{\bf 1. A "Bezout-like" condition for multiples of} $C$ {\bf and of
lines.}

\medskip
\medskip
{}From now on we assume that $I$ is an ideal of type (*) in $k[x_0,...,x_3]$,
i.e. that the points $P_i$ are in $\PP 3$.  Let $L$ be the line $P_iP_j$,
and let $n, t$
be  natural numbers; Proposition 1.1 will give a Bezout-type condition that
forces the elements of the linear system $I_t$ to contain the line $L$ with
multiplicity at least $n$ (it is actually just Bezout for $n = 1$).\par
Assuming
further that the $P_i$'s lie on a rational normal curve $C$,
Proposition 1.3 gives
an analogous condition that forces the elements of $I_t$ to contain the
curve $C$ with multiplicity at least $n$.
\par
Let  $(x)^+ = max\{ x,0\}$. We have
\par
\medskip

\noindent {\bf Proposition 1.1:} {\it Let} $I =
{\goth p}^{m_1} _1 \cap \ldots
\cap {\goth p}^{m_s} _s\  $ {\it be an ideal of type (*), and }$L$ {\it be
the line} $P_iP_j$. {\it If }$I_L$ {\it is the ideal of L, and}
$n\leq(m_i+m_j-t)^+$, {\it  then} $I_t \subseteq (I_L ^n)_t$.
\par
\medskip
{\it Proof:} The statement is obvious for $n=0$ and $n=1$, so let $n>1$.
 We may assume
$P_i=(0:0:0:1), P_j=(0:0:1:0)$, hence $I_L=(x_0, x_1)$. Let $f\in I_t$: by
 Bezout's theorem
applied to the intersection of $\{ f=0\}$ with the plane  $\{ x_1-ax_0=0\}$,
$a \in k$,
it is easy to prove that $f(x_0,ax_0,x_2,x_3) \in (x_0)^n, \forall a \in k$.
Hence $f(x_0,x_1,x_2,x_3) \in (x_0,x_1)^n$.\par
\qquad \qquad \qquad \qquad \qquad \qquad
\qquad \qquad \qquad \qquad \qquad \qquad
\qquad \qquad \qquad \qquad \qed
\par
\medskip
\medskip
\noindent {\bf Definition 1.2:} {\it Let} $I$ {\it be as above, and } $t,n \in
{\Bbb N}$. {\it We will say that} $I_t$ {\it satisfies property} ${\cal P}(n)$
{\it if and only if} $\quad \forall l,\quad 1\leq l\leq n:\quad 3t+5(1-l) <
\sum
_{i=1}^s (m_i -l+1)^+$.
\par
\medskip
\medskip
\noindent {\bf Proposition 1.3:} {\it Let} $I ={\goth p}^{m_1} _1 \cap \ldots
\cap {\goth p}^{m_s} _s\  $ {\it be an ideal of type (*)  such that the}
$P_i${\it 's are distinct points on a rational cubic curve C  and let} $t,n \in
{\Bbb N}$ {\it be such that} $I_t$ {\it satisfies property} ${\cal P}(n)$. {\it
Then, if} $I_C$ {\it is the ideal of C, we have that} $I_t \subseteq
(I_C^n)_t$.
\par
\medskip
{\it Proof:} Note that $I_C ^n = I_C ^{(n)}$, i.e. $I_C ^n$ is saturated and
represents the
$n^{th}$ infinitesimal neighborhood of $C$. In fact, by [{\bf 5}], Cor. 2.2,
we have:
\par
\medskip

\noindent {\bf Corollary} (Robbiano): {\it Let I be the ideal associated to a
complete intersection of codimension} $\leq 2$  {\it inside the Segre
embedding of} $\PP 1\times \PP {r-1}$  {\it in} $\PP {2r-1}$.
\par
{\it Then} $I^n$ {\it is primary for every n.}
\par
\medskip
{\it } Since the ideal of a rational normal cubic $C \subseteq \PP 3$ can be
obtained
by intersecting the ideal of $\PP 1\times \PP 2 \subset \PP 5$ which is
given by the maximal minors of a matrix of type $\pmatrix
{ x_0& x_1 & x_2\cr x_5& x_4& x_3}$, with the hyperplanes $\{ x_1 = x_5\}$,
$\{ x_2=x_4\}$, we get what we want (see also [{\bf A-S-V}], 6.9).\quad \qed
\par
\medskip
Let $X$  be the blow-up of $\PP 3$ along $C$. Then we have  ${\rm Pic} X
\cong {\bf Z}\oplus {\bf Z}$, and we can choose as generators the
exceptional divisor $E$ and the divisor $H$, corresponding to the strict
transform of a generic plane of  $\PP 3$.
\par
\medskip
Let $S$ be a surface in the linear system $I_t$ and let $S^\prime \subset
X$ be its strict transform.  Then Proposition 1.3 is equivalent to:
\par
\medskip

\noindent {\bf Proposition 1.4:} {\it Let} $I,n,t$ {\it be as in Proposition
1.3. Then} $S^\prime \in \vert tH - nE\vert$.
\medskip
\medskip
{\it } It is well known that $E \cong$ ${\Bbb P}({\cal N}_C)$
(e.g. see [{\bf 4}]), hence, since $ {\cal N}_C \cong {\cal O}_C(5)\oplus
{\cal O}_C(5)$, we have that
$$ E \cong {\Bbb P}({\cal N}_C) \cong {\Bbb P}({\cal O}_C(5) \oplus
{\cal O}_C(5)) \cong {\Bbb P}({\cal O}_C\oplus {\cal O}_C)$$
\noindent is isomorphic to a quadric surface
${\Bbb P}({\cal O}_{\PP 1}\oplus {\cal O}_{\PP 1})$.
\par
Let $\pi : X \rightarrow \PP 3$ be the canonical projection; then one of the
two rulings of $E$ is given by the lines $L_P = \pi ^{-1}(P)$, $P \in C$, and
the other is given by the zero-loci $C^\prime$ of sections of ${\cal N}_C$
($\pi (C^\prime) = C$ ).
\par
Let us use the notation $(a,b)$ for the divisor class of $aL_P + bC^\prime$;
then we have that $H\cdot E $, as divisor on $E$, is
$(E\cdot H)\vert _E =(3,0)$.
\par
Of course $H\cdot H = H^2$ is the strict transform of a generic line of
 $\PP 3$ (not touching $C$), so $H^2\cdot E = 0$.
\par
In order to determine the  $E^2$, consider $E^2\cdot H$. Since $E^2 =
E\vert _E$, we have $ E^2\cdot H = E\cdot (E\cdot H) = (E\cdot H)\vert _E
\cdot E^2 $. Let $E^2 = (a,b)$.  Since we have also $E^2\cdot H = E\vert _H
\cdot E\vert _H $, and $H \cong$ $\{$ the blow-up of $\PP 2$ at three
points $\}$, where  $E\vert _H$ is the exceptional divisor of such a blow
up, we get that $E^2 \cdot H = (E\vert _H \cdot E\vert _H) = -3$. Hence
$(H\cdot E)\vert _E \cdot E^2 = (3,0)\cdot (a,b) = 3b$,
which implies \  $b = -1$.
\par
In order to determine $a$, consider instead: $ E^3 = E^2\vert _E = E\vert _E
\cdot E\vert _E = -2a$ (since $E\vert _E = (a,-1)$ and $(a,-1)\cdot (a,-1)
= -2a$).
\par
On the other hand, $(3H-E)^3 = 27H^3 - 27H^2\cdot E + 9H \cdot E^2 - E^3$,
and so, since $H^2\cdot E = 0$, $H^3 = 1$ and $H\cdot E^2 = -3$, we have
$(3H-E)^3 = -E^3$.
\par
It is not hard to compute $(3H-E)^3$: this is the number of intersections of
the strict transforms $S_1 ^\prime$,$S_2 ^\prime$, $S_3 ^\prime $ of
three generic cubic surfaces $S_1$,$S_2$,$S_3$ containing $C$ in $\PP 3$.
\par
Consider $S_1 ^\prime \in \vert 3H-E\vert$.  The cubic surface $S_1$ in
$\PP 3$ is isomorphic to the blow-up of $\PP 2$ in six generic points
$P_1,...,P_6$, so Pic$S_1 \cong {\bf Z}^7$, and we can write
$(t;m_1,...,m_6)$ for the divisor class in $S_1$ of the strict transform of
a curve of degree t in $\PP 2$ which has multiplicity $m_i$ in $P_i$,
$i=1,...,6$; e.g. we have that a plane section of  $S_1$ is in
$(3;1,1,1,1,1,1)$.
\par
With this notation, we can assume  $C \in (1;0,0,0,0,0,0)$, and we have
$S_2\cdot S_1 \in (9,3,3,3,3,3,3) $ on $S_1$, with $S_2\cdot S_1 = C \cup
C^\prime$ and $C^\prime \in (8;3,3,3,3,3,3)$ ($C^\prime$ is a sextic curve
of genus 3). So the number we want is $S_3\cdot C^\prime - C\cdot
C^\prime = 18 - 8 = 10$.
\par
Hence $E^3 = -10$ \   and \   $a = 5$.
\medskip
\medskip
Now we want to prove Proposition 1.4.
Let $Z = (P_1,...,P_s;m_1,...,m_s)$ be the scheme defined by $I$ (hence all
$P_i \in C$), and $S \subset \PP 3$ be a surface of degree $t$ which
contains $Z$.
\par
\noindent  Note that for $n=1$ the statement is trivial by Bezout Theorem.
\par
Let us work by induction on $n$, the case $n=1$ being done. Suppose the
proposition is true for $n-1>0$, and let us check that it is true for $n$.
\par
Since $3t+5(1-l) < \sum _i ^s (m_i -l+1)^+$ for $1\leq l \leq n-1$, then, by
induction, any surface of degree $t$ containing $Z$ is such that its strict
transform is in $\vert tH-(n-1)E\vert$. \par
\noindent So, let $S^\prime $ be the strict transform of $S$ in $X$; then
$S^\prime \in \vert tH-(n-1)E\vert$. Consider (on $E$)
$$S^\prime \cdot E = (tH-(n-1)E)\cdot E = (3t,0) - (n-1)(5,-1) = (3t-5(n-1),
(n-1)).$$
Let $L_i = \pi ^{-1}(P_i)$; then $\pi ^{-1}S = S^\prime + (n-1)E$ has to
contain the divisors $m_iL_i \in (m_i,0)$, $i=1,...,s$ of $E$, hence
 $S^\prime $ has to contain $(m_i- n + 1)^+L_i$, i.e. at least a divisor
$(\sum _i ^s (m_i-n+1)^+,0)$ on $E$.
\par
So, if $3t-5(n-1) < \sum _i ^s (m_i-n+1)^+$, $E$ has to be a fixed
component in $tH-(n-1)E$, i.e. $S^\prime \in \vert tH-nE\vert$.
\qquad \qed
\par
\medskip
It can be of some interest to give an actual computation of
$\dim (I_C ^n)_t$:
\medskip
\medskip

\noindent {\bf Proposition 1.5:} {\it Let } $I_C \subseteq k[x_0,...,x_3]$ {\it
be the ideal of a rational normal curve} $C \subseteq \PP 3$. {\it Then we
have:}  \par $\dim (I_C^n)_t = 0${\it , for} \ \  $t \leq 2n-1${\it , and} \par
$\dim (I_C^n)_t = {t+3 \choose 3} - {n+1 \choose 2}(3t+6) +
5(1+...+n^2)${\it , \  for} \ $t \geq 2n-1$; {\it moreover in this case, if}
${\cal I}_C$ {\it is the ideal sheaf associated to} $I_C${\it , we have :}
$H^1(\PP 3,{\cal I}_C^n(t)) = 0$.
\par
\medskip
\medskip
{\it Proof:} The case $t \leq 2n-1$ is obvious
since $I_C ^n $ is generated in degree $2n$, so consider $t \geq 2n$.
\par
Let $q=0$ be the equation of a smooth quadric $Q$ containing
$C$; then multiplication by $q$ defines an injection
$0 \rightarrow I_C^{n-1} \rightarrow I_C^n$, from which,  sheafifing, we get:
$$0 \rightarrow {\cal I}_ C ^{n-1}(-2) \rightarrow {\cal I}_C ^n \rightarrow
{\cal O}_Q(-nC) \rightarrow 0.$$
Since $C$ is of type $(2,1)$ as a divisor on $Q$, twisting  by
${\cal O}_{\PP 3}(t)$ we get
$$0 \rightarrow {\cal I}_ C ^{n-1}(t-2) \rightarrow {\cal I}_C ^n(t)
\rightarrow {\cal O}_Q(t-2n,t-n) \rightarrow 0$$
\noindent which, passing to cohomology, yields:\par
\medskip
\noindent $0 \rightarrow H^0(\PP 3,{\cal I}_ C ^{n-1}(t-2)) \rightarrow
H^0(\PP 3,{\cal I}_C ^n(t)) \rightarrow H^0(Q,{\cal O}_Q(t-2n,t-n))
\rightarrow $
\par
$\qquad \qquad \rightarrow H^1(\PP 3,{\cal I}_ C ^{n-1}(t-2))
\rightarrow H^1(\PP 3,{\cal I}_C ^n(t)) \rightarrow
H^1(Q,{\cal O}_Q(t-2n,t-n))  $
\par
\noindent where $H^1(Q,{\cal O}_Q(t-2n,t-n)) = 0$ since $t \geq 2n$.
\par
\medskip
Now we work by induction on $n$.  For $n=1$, it is well known that
\par
$\qquad \quad \dim H^0(\PP 3,{\cal I}_C ^n(t)) = {t+3 \choose 3} -
(3t+1)\quad $  and $\quad H^1(\PP 3,{\cal I}_C (t)) = 0$.\par
So, suppose $n \geq 2$. Since $t-2 \geq 2(n-1)$, we have
$H^1(\PP 3,{\cal I}_C ^{n-1}(t-2)) = 0$ by induction hypothesis, hence
\noindent $H^0(\PP 3,{\cal I}_C ^n(t)) = H^0(\PP 3,{\cal I}_C ^{n-1}(t-2)) +
H^0(Q,{\cal O}_Q(t-2n,t-n)) =$ \par
$= {t+1 \choose 3} - {n \choose 2}(3(t-2)+6) + 5(1+...+(n-1)^2) +
(t-n+1)(t-2n+1) =$
\par
$\qquad \qquad = {t+3 \choose 3} - {n+1 \choose 2}(3t+6) + 5(1+...+n^2).$
 \quad \qed \par
\medskip
We conclude this section with the following result, which gives, in case
$n=1$, a sort of inverse with respect to Proposition 1.3.
\par
\medskip

\noindent {\bf Proposition 1.6:} {\it Let} $I$ {\it be as in Proposition 1.3,
and let} $I_t \neq \{0\}$. {\it Then} $I_t \subseteq (I_C)_t$
{\it if and only if } $3t < \sum _i ^s m_i$.
\par
\medskip
{\it Proof:} By proposition 1.3, it remains only to prove that
$\{0\} \neq I_t \subseteq (I_C)_t$ implies $3t < \sum _{i=1} ^s m_i$.
This follows from the fact that the inequality
$3t \geq \sum _{i=1} ^s m_i$ allows to find a surface in $I_t$
made of planes (each of them passing at most through three
of the $P_i$'s),
hence not contained in $(I_C)_t$. \quad \qed
\medskip
\medskip
\noindent {\bf Remark 1.7:} Notice that it is not possible, instead, to do the
same with Proposition 1.1, in fact the condition $t \geq m_i+m_j$ does not
guarantee that the line $P_iP_j$ is not fixed for $I_t$. \par
For instance, let $Z=(P_1,...,P_7;3,3,2,2,2,2,1)$ with
$P_1,...,P_7$ on $C$; we have that
 the only surface of degree 4 containing $Z$ is given by the union of the
two quadric cones with vertices in $P_1,P_2$ which contain $C$.
\par
Hence the lines $P_1P_7$ and $P_2P_7$ are fixed for $I_4$, even if
$4=t=m_1+m_7=m_2+m_7$.
\par
\medskip
On the other hand, the two numerical conditions ($t \geq m_i+m_j$ and
$t \geq \sum _i ^s m_i$) together are equivalent to the two geometric
conditions (see Corollary 2.3).
\par

\medskip
\medskip
\medskip

\noindent {\bf 2. The algorithm to compute the Hilbert function of fat points
on a cubic curve.}

\medskip
Let $I$,$J$ be, respectively, the homogeneous ideals of the schemes
$Z$,$W \in \PP 3$ of fat points $$Z = (P_1,...,P_s;m_1,...,m_s),\quad
W = (P_1,...,P_s;m_1,...,m_{s-1},m_s +1)$$ where, from now on,
$P_1,...,P_s$
are on rational normal cubic $C$, and \par
\centerline{$m_1 \geq m_2 \geq \ldots \geq m_s \geq 0$.}
\par
We want to give a method that can compute $\dim J_t$ (for every
$t \geq 0$) from the data:  $t,m_1,...,m_s$ and $\dim I_t$.
The result will not depend on the position of the points on $C$.
\par
This will answer to our question, since one will be able to compute
$\dim J_t$ working by recursion on $s$ and $m_s$.
The algorithm will be given for $s\geq 2$, since the case $s=1$
is quite trivial. In fact for $s=1$ we have $Z = (P,m)$, $W = (P,m+1)$,
and so we get that $t \geq m$ implies:
$$\dim (I/J)_t = \dim I_t - \dim J_t =
{m+3 \choose 3} - {m+2 \choose 3} = {m+2 \choose 2};$$
while when $t < m$ , trivially $\dim (I/J)_t = 0$.
\par
Hence, in the sequel, we will always suppose $s\geq 2$.
\medskip
We will determine  $\dim (I/J)_t$ via  a scheme of fat points $N \subset \PP
2$.
\medskip

\noindent {\bf Definition 2.1.} {\it Let} $Z = (P_1,...,P_s;m_1,...,m_s)$
{\it be a scheme of fat points in} $\PP 3$ {\it with support on a rational
normal curve C and let I be the corresponding homogeneous ideal.
We say that } $N = (Q,Q_1,\ldots ,Q_{s-1};n,n_1,
\ldots ,n_{s-1}) \subseteq \Pi \cong \PP 2$ {\it is the t-projection of}
$Z $ {\it from} $P_s${\it , if the points} $Q_1,\ldots ,Q_{s-1}$
{\it are the projection from } $P_s$ {\it of } $P_1,\ldots ,P_{s-1}$
{\it on a plane} $\Pi$ {\it not containing} $P_s${\it , while} $Q$
{\it is the projection of } $P_s$ {\it itself along the tangent line to}
$C$ {\it at} $P_s$ {\it and the numbers } $n,n_1,...,n_{s-1}$
{\it are defined as follows:}
$$n_i = (m_i + m_s - t)^+ , \qquad
n={\rm min}\{m_s+1,{\rm sup}\{ \nu \in {\Bbb N}\vert {\cal P}(\nu )\
{\rm holds\  for}\  I_t\}\}.$$
\medskip
We can always suppose $P_s = (0:0:0:1)$, and $\Pi $ to be $\{x_3=0\}$.
\par
Of course the points  $Q,Q_1,\ldots ,Q_{s-1}$ lie on the conic
$\Gamma $ which is the projection of $C$ from $P_s$,
so the Hilbert function of $N$ is known (see  [{\bf 2}]).
\par
Our result is:\par
\medskip

\noindent {\bf  Theorem 2.2:} {\it Let} $I$, $J$ {\it be respectively
the homogeneous ideals of the schemes of fat points}
$$Z = (P_1,...,P_s;m_1,...,m_s),\quad
W=(P_1,...,P_s;m_1,...,m_{s-1},m_s +1)$$
{\it where the} $P_i${\it 's are distinct points of} $\PP 3$ {\it lying
on a rational normal cubic} $C$, {\it and} $m_1 \geq m_2 \geq \ldots
\geq m_s \geq 0$, $s\geq 2$.
{\it We have, for every} $t\geq 0${\it , that} $ \dim (I / J)_t$
{\it equals the dimension, in degree} $m_s${\it , of the ideal}
$I_N \subseteq k[x_0,x_1,x_2]$ {\it of the t-projection}
$N$ {\it of } $Z$  {\it from} $P_s$, {\it i.e. :}
$$ \dim (I/ J)_t = \dim (I_N)_{m_s} .$$
\medskip
Note that $L_i = \{{\rm line} \  P_iP_s\}$ and $C$ are
fixed multiple curves for the surfaces in the linear system $I_t$
with multiplicities at least $n_i, n$, respectively (see Proposition 1.1
and 1.3).
The theorem shows their role in determining the Hilbert Function
of $W$; in particular, when
$n=n_i=0$, i.e. when $N = \emptyset $, the difference between
$\dim I_t$ and $\dim J_t$ is ${m_s + 2 \choose 2}
= {m_s + 3 \choose 3} - {m_s + 2 \choose 3}$,
i.e. it is what it "should be", in the sense that passing to
multiplicity $m_s + 1$ on $P_s$ imposes exactly
${m_s + 2 \choose 2} $ new independent conditions to
surfaces of degree $t$.
\par
Thus, we have:
\par
\medskip

\noindent {\bf Corollary 2.3:} {\it For any ideal} $I$ {\it as in
Theorem 2.2, if} $m_s>0$ {\it , then the following are equivalent:}
\par
\noindent {\it i) } $I_t$ {\it is regular, (i.e. the fat points impose
independent conditions to surfaces of degree t);}
\par
\noindent {\it ii) neither} $C$ {\it nor,} $\forall i\in \{1,..,s\},${\it any
of the lines} $P_iP_j$, $j\neq i,$ {\it is a fixed locus for}
$(I:{\goth p}_i)_t = ({\goth p}^{m_1} _1 \cap
\ldots \cap {\goth p}^{m_i-1} _i
\cap \ldots \cap {\goth p}^{m_s} _s)_t $;
\par
\noindent {\it iii)}
$3t \geq \sum _{i=1}^s m_i-1$  {\it and}  $t \geq m_1+m_2-1$.
\par
\medskip
\noindent {\it Proof:} By Bezout's Theorem, {\it ii)} implies {\it iii)}.
Let us show now that {\it iii)} implies the regularity of $I_t$.
Consider that one can get $I_t$ starting from $({\goth p}_1^{m_1})_t$,
which is regular, and "adding the multiplicities on the $P_i$'s
one at a time", i.e. considering the ideals associated to the schemes
$$(P_1,...,P_s;m_1,0,...,0),\quad (P_1,...,P_s;m_1,1,0,...,0),\quad
\quad (P_1,...,P_s;m_1,2,...,0)$$
and so on.  At any step we have that the t-projection of such schemes
from the "last" point is empty, so Theorem 2.2  tells us that adding one
to the multiplicity of the last point imposes  independent conditions,
and the system remains regular. (See also [{\bf 3}]).\par
In order to prove that {\it i)} implies {\it ii)}, suppose that either
$C$ or a line $P_iP_j$ are fixed components for $(I:{\goth p}_i)_t$
for some $i=1,...,s$. Then we will show that
$\dim ((I:{\goth p}_i)/I)_t < e = {m_i+2\choose 2}$,
hence that $I_t$ cannot be regular.
\par
If $I_t$ were regular, let $P_i=(0:0:0:1)$; then there would exist
$F_1,...,F_e \in (I:{\goth p}_i)_t$ such that, locally,
they generate $(x,y,z)^{m_i}$ modulo $I_t$, hence they are of type
(in affine coordinates):
$$F_1 = \tilde F_1 + x^{m_i};
\quad F_2 = \tilde F_2 + x^{m_i-1}y;
\quad ... ;\quad F_e = \tilde F_e + z^{m_i}$$
\noindent where the $\tilde F_i $ have degree $\geq m_i+1$.
Now, if $C = \{x=t,y=t^2,z=t^3\}$ is fixed for $(I:{\goth p}_i)_t $,
the above equations should become identities in $t$, but this is clearly
impossible for the first one, since $x^{m_i} = t^{m_i}$ while
$\tilde F_1$ has degree $\geq m_i+1$ in $t$.
\par
If the fixed component is a line $P_iP_j$, we work in the same way,
e.g. assuming that the line is given by $\{x=t,y=0,z=0\}$.\par
\qquad \qquad \qquad \qquad \qquad \qquad
\qquad \qquad \qquad \qquad \qquad
\qed
\medskip
\medskip
\medskip

\noindent {\bf 3. Preliminary Lemmata.}

\medskip
The proof of Theorem 2.2 will be given showing first that
$\dim (I / J)_t \leq \dim (I_N)_{m_s}$ (Lemma 3.1).
Then we will consider several
particular cases, with which we will deal with lemmata 3.2 to 3.5.
This will leave us only with cases in Remark 3.6. \par
Lemmata 3.7 to 3.10
describe geometric properties of the cases listed by the Remark, and
they will be used in the next section for the proof of the theorem.
\medskip
 From now on, we will always suppose that $I$, $J$, $I_N$ are as in Theorem 2.2
 and  $s \geq 2$.
\medskip

\noindent {\bf Lemma 3.1:} {\it For every} \ $t \geq 0$,
{\it we have}\  $\dim (I / J)_t \leq \dim (I_N)_{m_s}.$
\medskip
\noindent {\it Proof:} Let $P_s = (0:0:0:1)$.
If $F = F_{m_s} (x_0,x_1,x_2) x^{t-m_s}_3 + F_{m_s+1}
(x_0,x_1,x_2)x_3^{t-m_s-1}+ \ldots  +F_t (x_0,x_1,x_2)$
is a form of $I_t$, then it is easy to prove that $F_{m_s}$,
i.e. a polynomial defining the tangent cone to $\{ F=0 \}$ in $P_s$,
is a form of $(I_N)_{m_s}$. Since the application $(I/J)_t \rightarrow
(I_N)_{m_s}$ that maps the class of $F$ to $F_{m_s}$ is injective,
we have the conclusion.\par
\qquad \qquad \qquad \qquad \qquad
\qquad \qquad \qquad \qquad \qquad
 \qquad \qquad \qquad \qquad \qquad \qquad
\qed
\medskip
\medskip

\noindent {\bf Lemma 3.2:} {\it If} $n = m_s+1$ {\it or }
$n_1 \geq m_s+1$,  {\it then} $\dim (I / J)_t =
 \dim (I_N)_{m_s} = 0.$
\medskip
\noindent {\it Proof:} Trivially $\dim (I_N)_{m_s}= 0$,
 so, by Lemma 3.1, we are done.
 \qquad  \qquad
\qed \par
\medskip
\medskip

\noindent {\bf Lemma 3.3:} {\it If}\  $n = 0$ {\it and} $n_1 = 0$,
{\it then }
$\dim (I/J)_t = \dim (I_N)_{m_s} = {m_s + 2 \choose 2}$
\par
\medskip
\noindent {\it Proof:}  If $m_{s-1} \geq m_s+1$,
the result follows from \ref{3}, Proposition 5.
 For $m_{s-1} = m_s$ it is easy to extend the
above proposition to our case.
\qquad
\qed
\medskip

\noindent {\bf Lemma 3.4:} {\it Let}  $n = 0$, $n_1> 0$ {\it and one of the
following cases occurs:}\par
\noindent $a) \quad s = 2$,\quad {\it and} $\quad n_1 \leq m_2$.
\par
\noindent $b) \quad s = 3$ {\it or} $4$,\quad {\it and} $\quad n_1+n_2 \leq
m_s$.
\par
 {\it Then} \quad $\dim (I / J)_t = \dim (I_N)_{m_s}$.
\medskip
\noindent {\it Proof:} Note that in both cases $\dim (I_N)_{m_s}$ is
known since  $(I_N)_{m_s}$ is regular
(e.g. see [{\bf 2}]). \par
In case {\it a)} suppose $P_1 = [0:0:1:1]$ and
$P_2 = [0:0:0:1]$, so $I_t  \subseteq (x,y)^{n_1}$. \par
We have to show that $\dim (I / J)_t \geq  \dim (I_N)_{m_s}$,
since we have seen (Lemma 3.1)  that the opposite inclusion always holds.
\par
Let us consider the monomials (in affine coordinates): \par
\medskip
\noindent $x^{m_2}, x^{m_2-1}y, x^{m_2-2}y^2,...,y^{m_2}$;\par
\noindent $zx^{m_2-1}, zx^{m_2-2}y,...,zy^{m_2-1}$;\par
$\quad \quad \ldots $\par
\noindent $z^{m_2-n_1}x^{n_1}, z^{m_2-n_1}x^{n_1-1}y,...,
z^{m_2-n_1}xy^{n_1-1}, z^{m_2-n_1}y^{n_1}$.\par
\medskip
All these monomials have multiplicity $m_2$ at $P_2$, and
at least $n_1$ at $P_1$.  Since $t-m_2 = m_1-n_1$, multiplying
the above monomials by $(z-1)^{t-m_2}$ we get linearly independent polynomials
of degree $t$ with multiplicity at least $m_1$ at $P_1$ and exactly
$m_2$ at $P_2$ (where they have independent initial forms), hence:
$$\dim (I / J)_t \geq (m_2+1) + m_2 + (m_2-1)
+ ... + (n_1+1) = {m_2+2 \choose 2} - {n_1+1 \choose 2} = \dim (I_N)_{m_2}.
$$
So we are done in case {\it a)}.
\par
In case $b)$, for $s=3$,  let
$P_1 = [0:0:1:1]$, $P_2 = [1:0:0:1]$ and $P_3 = [0:0:0:1]$.
\par
If $n_2 = 0$, we consider the same monomial as above,
but with $m_3$ instead of $m_2$; multiplying them
by $(x-z-1)^{t-m_3}$ we get
$\dim (I / J)_t = \dim (I_N)_{m_s}$, as we did above.
 \par
If $n_2 > 0$, we consider instead the monomials:
\par
\medskip
\noindent $x^{m_3-n_2}y^{n_2},..., xy^{m_3-1}, y^{m_3}$;\par
\noindent $zx^{m_3-n_2}y^{n_2-1},..., zy^{m_3-1}$;\par
$\quad \quad \ldots $\par
\noindent $z^{n_2}x^{m_3-n_2}y^0,..., z^{n_2}y^{m_3-n_2}$;\par
\noindent $z^{n_2+1}x^{m_3-n_2-1},..., z^{n_2+1}y^{m_3-n_2-1}$;\par
$\quad \quad \ldots $\par
\noindent $z^{m_3-n_1}x^{n_1},..., z^{m_3-n_1}y^{n_1}$.\par
\medskip
With the same kind of reasoning as before, we get \par
\medskip
$\qquad \qquad
\dim (I / J)_t \geq (m_3 - n_2 +1)n_2 +
{(m_3-n_2+1+n_1+1)(m_3-n_2-n_1+1) \over 2} = $
\par
$\qquad \qquad \qquad \qquad ={m_3+2\choose 2}-
{n_2+1 \choose 2} - {n_1+1 \choose 2} = \dim (I_N)_{m_3}.
$ \par
\medskip
\noindent The case $s=4$ is completely analogous by taking $P_1 = [0:0:1:1]$,
$P_2 = [1:0:0:1]$, $P_3 = [0:1:0:1]$ and $P_4 = [0:0:0:1]$,
so we leave it to the reader.
\par
\qquad \qquad \qquad \qquad \qquad \qquad \qquad \qquad
\qquad \qquad \qquad \qquad \qquad \qquad \qquad \qquad
\qquad \qed  \par
\medskip
\medskip

\noindent {\bf Lemma 3.5:} {\it Let } $n \leq 1$, $n_1 \leq 1$, $s \geq 3$,
$m_1=m_2=...=m_s=1$, {\it , then}
$\dim (I / J)_t = \dim (I_N)_{m_s}$.
\par
\medskip
\noindent {\it Proof:} If $3t <s$, we have $n = 1$ and ${\cal P}(2)$ does
 not hold for $I_t$. It follows that $3t-5 \geq 0$, so $t \geq 2$, $s \geq 7$,
$N=(Q;1)$. Hence $\dim I_t = 3$, $\dim J_t = 1$, $\dim (I_N)_1 = 2$
and we are done. \par
Let $3t\geq s$, so $n = 0$. For $t = 1$ we have $s = 3$,
$N=(Q_1,Q_2;1,1)$ so $\dim I_1 = 1$, $\dim J_1 = 0$,
and $\dim (I_N)_1 = 1$. For $t >1$,  we have $n_1=0$,
and the conclusion follows by Lemma 3.3.\par
\qquad  \qquad \qquad \qquad \qquad \qquad \qquad
\qquad \qquad \qquad \qquad \qquad \qquad \qed
\medskip

 In the following remark we list the cases not covered by the previous
 lemmata  (notice that for $2 \leq s \leq 5$
we only have $n = 0$ or $n = m_s+1$):
\par
\medskip
\medskip

\noindent {\bf Remark 3.6:} {\it In order to complete the proof of
Theorem 2.2 the following cases remain to be considered, where we always
have}  $\quad s \geq 3$, $\quad m_1 \geq 2$,
 $\quad n \leq m_s$,$\quad n_1 \leq m_s$,$\quad n+n_1>0$
{\it (hence} $t \geq m_1$ {\it and } $m_s>0${\it ):}
\par
\medskip
\noindent ${\bf 1)} \quad n_1+n_2 \geq m_s+1$;\par
\noindent ${\bf 2)} \quad s\geq 6, \quad n+n_1\geq m_s+1, \quad m_1 > m_s$;\par
\noindent ${\bf 3)} \quad s\geq 6, \quad  n+n_1\geq m_s+1, \quad m_1 = m_s>1$;
\par
\noindent ${\bf 4)} \quad s \geq 5 , \quad  m_s \geq n+n_1 ,
\quad  m_s \geq n_1+n_2$.\par

\medskip
\medskip
\medskip

\noindent {\bf Lemma 3.7:} {\it In case 1) we have:}\par
\noindent {\it a) \quad The line } $Q_1Q_2$ {\it is a fixed component for}
$(I_N)_{m_s}$;
\par
\noindent {\it b) \quad The plane } $P_1P_2P_s$
{\it is a fixed component for} $I_t $, $J_t$.
\par
\medskip
\noindent {\it Proof:} Point $a)$ is obvious. For $b)$,
since $n_2\geq 1$, notice that the surfaces of $I_t $ contain the line $P_1P_2$
with multiplicity $m_1+m_2-t$,  the line $P_1P_s$  with multiplicity
$m_1+m_s-t=n_1$ and the line $P_2P_s$
 with multiplicity $m_2+m_s-t=n_2$, hence the degree of intersection
of the plane $P_1P_2P_s$ with those surfaces is: \par
\noindent $m_1+m_2-t+n_1+n_2 =
n_1+n_2 - 2m_s + 2t - t + n_1 + n_2 = 2(n_1+n_2-m_s)+t \geq t+2$
\par
\noindent so, by Bezout, the plane has to be a fixed component. \quad \qed
\par
\medskip
\medskip

\noindent {\bf Lemma 3.8:} {\it In cases 2) and 3) we have:}\par
\noindent {\it a) \quad The line } $Q_1Q$ {\it is
a fixed component for} $(I_N)_{m_s}$;
\par
\noindent {\it b) \quad The quadric cone on} $\Gamma$
{\it with vertex in } $P_1$ {\it is
a fixed component for} $I_t $, $J_t$.
\par
\medskip
\noindent {\it Proof:} As before, $a)$ is trivial, while $b)$
 follows by Bezout, considering the fact that the surfaces in
 $I_t $ contain the curve $C$ with multiplicity at least $n$
 and the lines $P_1P_i$, $i=2,...,s,$ with multiplicity
 at least $m_1+m_i-t$. Since $n_1 \geq 1$, and ${\cal P} (n)$
holds for $I_t,$
one gets that their multiplicity of intersection with the cone is
$\geq 3n+(s-2)m_1+{\sum _{i=1} ^s m_i}-(s-1)t >
3n+(s-2)m_1+(s-5)(n-1) +3t-(s-1)t=
(s-2)(n+m_1-t-1)+2t+3=(s-2)(n+n_1-m_s-1)+2t+3 \geq 2t+3$.
 \quad \qed \par
 \medskip
\medskip

\noindent {\bf Lemma 3.9:} {\it In case 3) we have:}\par
\noindent {\it a) \quad The conic } $\Gamma $ {\it is
a fixed component for} $(I_N)_{m_s}$;
\par
\noindent {\it b) \quad The quadric cone on} $\Gamma$
{\it with vertex in } $P_s$ {\it is
a fixed component for} $I_t $, $J_t$.
\par
\medskip
\noindent {\it Proof:} Since $m_1=m_s$, all the cones on $\Gamma $
with vertex in a $P_i$ are fixed for $I_t $, by the previous Lemma, so $b)$
is done. \par
To show $a)$ we will check, if $m_1=...=m_s=m$, that  \par
\noindent $(1) \qquad \qquad \qquad \qquad \qquad (s-1)n_1+n-1-2m \geq 0$.\par
Since $n_1 \geq 1$ and ${\cal P} (n)$ holds for $I_t,$ we have:\par
\noindent $(s-5)n < sm-3t+s-5 = sm-3(2m-n_1)+s-5 \leq
(s-5)(n+n_1-1)+3n_1+s-5-m = (s-5)n+(s-2)n_1-m\ $ thus, \par
\noindent $(2) \qquad \qquad \qquad \qquad \qquad \qquad m < (s-2)n_1$. \par
On the other hand, since ${\cal P} (n+1)$ doesn't hold for $I_t,$
we have:\par
\noindent $  (3) \qquad \qquad \qquad \qquad (s-5)n \geq sm -3t =
sm-3(2m-n_1)$.
\par
To prove $(1) $, let us multiply it by $(s-5)$. We get: \par
\noindent $(s-5)((s-1)n_1+n-1-2m) \geq ^{{\rm (by} (3) {\rm )}}\
sm-3(2m-n_1)+(s^2-6s+5)n_1-s+5-2m(s-5) =
 (s-4)((s-2)n_1-m-1)+1 \geq ^{{\rm (by} (2) {\rm )}}\  (s-4)(m+1-m-1)+1 = 1$.
\quad \qed \par
\medskip
\medskip

\noindent {\bf Lemma 3.10:} {\it In case 4) we have}
$2m_s \geq  \sum _{i=1} ^{s-1} n_i +n$. \par
 \medskip
\noindent {\it Proof:} For $s = 5$, we have $n=0$ and $2m_s \geq 2(n_1+n_2)
\geq \sum _{i=1} ^{s-1} n_i$. \par
Let $s \geq 6$. If $n_4=0$ the result is obvious  since
$2m_s \geq n + n_1 + n_1 + n_2 \geq  n + n_1 + n_2 + n_3$. \par
Now let $n_4>0$. The case $n=0$ is not possible; in fact since
$t=m_1+m_s-n_1= m_2+m_s-n_2 = m_3+m_s-n_3$ and $s \geq 6$, we get
$3t = 3m_s+m_1+m_2+m_3-n_1-n_2-n_3 < \sum ^s _{i=1} m_i$. \par
So, $n_4>0$, $n>0$. Let $r=max\{i \vert n_i>0\}$,
$4\leq r\leq s-1$; then  (since
$t = m_1+m_s-n_1 = ... = m_r+m_s-n_r$), we have : \par
\noindent $(4) \qquad \qquad \qquad \qquad
rt = \sum ^r _{i=1} m_i + rm_s - \sum ^r _{i=1} n_i$.
\par
Since ${\cal P} (n+1)$ doesn't hold for $I_t,$ we know that
$3t - \sum ^s _{i=1} m_i + (s-5)n \geq 0$, so: \par
\noindent $(s-5)n \geq \sum ^s _{i=1} m_i - 3t \geq
\sum ^r _{i=1} m_i + (s-r)m_s-3t =^{{\rm by (4)}} rt -rm_s +
\sum ^r _{i=1} n_i
+(s-r)m_s - 3t  =
(r-3)(m_1+m_s-n_1)+(s-2r)m_s+\sum ^r _{i=1} n_i \geq
(r-3)2m_s-(r-3)n_1+(s-2r)m_s+\sum ^r _{i=1} n_i  =
 (s-6)m_s-(r-3)n_1 + \sum ^r _{i=1} n_i \geq
(s-6)m_s-(s-4)n_1+\sum ^r _{i=1} n_i \geq
(s-4)(n+n_1)-(s-4)n_1+ \sum ^r _{i=1} n_i -2m_s.$
\par
{}From this we have: $\quad -n\geq \sum ^r _{i=1} n_i -2m_s$, hence the
conclusion.  \quad \quad \qed
\par
\medskip
\medskip
\medskip
\medskip

\noindent {\bf 4. Proof of Theorem 2.2.}

\medskip
The proof of the theorem works by induction on  $\sum ^s _{i=1} m_i$.
The first steps of the induction are covered by the lemmata 3.2 to
3.5 in \S 3; notice that also the trivial case $m_1=0$, i.e.
 $Z$ empty, is covered (by Lemma 3.3). \par
Now let us consider the cases left open in Remark 3.6.\par
\medskip
\medskip

\noindent {\bf Case 1) .}  Since $n_2 \geq 1$, we may consider
the following subscheme of $N$:
$$N^\prime = (Q,Q_1,...,Q_{s-1};n,n_1-1,n_2-1,n_3,...,n_{s-1}).$$
By Lemma 3.7, $(I_N)_t$ has the line $Q_1Q_2$ as a fixed component,
hence we have that $$\dim (I_N)_{m_s} = \dim (I_{N^\prime })_{m_s-1}.$$
\par
For a similar reason (the plane $P_1P_2P_s$ is a fixed component),
we have also:
$$\dim \left( {I \over J} \right) _{t} =
\dim \left( {I^\prime \over {J^\prime}} \right) _{t-1} $$
\par
\noindent
where $I^\prime$, $J^\prime$ are the ideals associated respectively
to the schemes
$$Z\prim =(P_1,...,P_s;m_1-1,m_2-1,m_3,...,m_{s-1},m_s-1),$$
$$W\prim = (P_1,...,P_s;m_1-1,m_2-1,m_3,...,m_{s-1},m_s).$$
Thus we are done, since it is quite easy to check that
$N^\prime$ is the $(t-1)$-projection from $P_s$ of the
scheme associated to $I^\prime $ and so,
by induction hypothesis, we have:
$$\dim \left( {I^\prime \over {J^\prime}} \right) _{t-1}
= \dim (I_{N^\prime} )_{m_s-1}.$$
 \qquad \qquad \qquad \qquad \qquad \qquad \qquad
\qquad \qquad \qquad \qquad \qquad \qquad \qquad \qquad \qed
\medskip

\noindent {\bf Case 2).} Note that $n_1 > 0$, $n>0$ .
In this case we proceed as above, but the fixed components we
"take away" are the quadric cone $\Lambda $ on  $\Gamma $ with vertex
in $P_1$ and the line $Q_1Q$ which are fixed for
$I _t$,$J_t$ and $(I_N)_{m_s}$ respectively (by Lemma 3.8).\par
So, let $I^\prime,J^\prime$ be the ideal associated to the schemes
$$Z^\prime = (P_1,...,P_s;m_1-2,m_2-1,...,m_{s-1}-1,m_s-1) =
(P_1,...,P_s;m_1{\prim},...,m_s{\prim}),$$
$$W^\prime = (P_1,...,P_s;m_1-2,m_2-1,...,m_{s-1}-1,m_s)$$ \par
\noindent respectively,
and let $$N\prim =(Q,Q_1,...,Q_{s-1};n-1,n_1-1,n_2,...,n_{s-1}).$$
\par
Then we have:
$$\dim \left( {I \over J} \right) _{t} =
\dim \left( {I^\prime \over {J^\prime}} \right) _{t-2} \qquad {\rm ;} \qquad
\dim (I_N)_{m_s} = \dim (I_{N^\prime })_{m_s-1}. $$ \par
In order to conclude, by induction hypothesis, we have to show that
$N\prim $ is the $(t-2)$-projection of $Z\prim $
from $P_s$, i.e. that the coefficients $n,n_1$ are the right
ones.
\par
In fact, $(m_1 -2) + (m_s -1) - (t-2) = n_1 - 1$; so, if we
let $n\prim $ be the
coefficient for $Q$ in the $(t-2)$-projection of $Z\prim $,
it only remains to check that $n-1=n\prim $.
\par
Since $\sum ^s _{i=1}m_i\prim = \sum ^s _{i=1}m_i - s - 1$, then for
$ 1\leq l \leq m_s-1$, we have
$3t+5(1-l)-\sum ^s _{i=1}(m_i-l+1)^+
= 3(t-2) + 5(1-(l-1)) - \sum ^s _{i=1}(m\prim _i - (l-1)+1)^+$.
It follows that $n \prim = n-1$.\qquad \qed
\medskip
\medskip

\noindent {\bf Case 3).} Here $m_1=...=m_s\geq 2$. By Lemma 3.9,
the conic $\Gamma $
is a fixed component for $(I_N)_{m_s}$, and the quadric cone
on $\Gamma $ with vertex in $P_s$ is fixed for $I_t,J_t$;
so, let $I^\prime,J^\prime $ be the ideal associated to the schemes
$$Z\prim=(P_1,...,P_s;m_1-1,m_2-1,...,m_{s-1}-1,m_s-2),$$
$$W\prim = (P_1,...,P_s;m_1-1,m_2-1,...,m_{s-1}-1,m_s-1)$$
 respectively,
and let $$N\prim = (Q,Q_1,...,Q_{s-1};n-1,n_1-1,...,n_{s-1}-1).$$ \par
We have $\dim (I/J)_t = \dim (I\prim /J\prim)_{t-2}$, and
$\dim (I_N)_{m_s} = \dim (I_{N\prim} )_{m_s-2}$, so if we show that
$N\prim $ is the $(t-2)$-projection of $Z\prim $ we
are done (by induction).\par
The kind of computations that are required are similar to the
ones shown in the previous cases, so we left them to the reader.
\qquad \qed
\par
\medskip

Note that in order to apply the inductive hypothesis to $W\prim$ and $Z\prim$
we had to use Lemma 3.9, and not Lemma 3.8 (which also
applies to this case).
\par
\medskip
In order to deal with the remaining case, we will show the following
lemmata:
\medskip
\medskip

\noindent {\bf Lemma 4.1:} {\it Let}
$\{h=0\}$ {\it be the plane}
$P_1P_2P_s$ {\it and let} $I^\prime,J^\prime$ {\it be the
homogeneous ideals associated respectively to the schemes:} \par
\noindent
$$Z\prim = (P_1,...,P_s;m_1-1,m_2-1,m_3,...,m_{s-1},m_s-1),$$
$$W\prim =(P_1,...,P_s;m_1-1,m_2-1,m_3,...,m_{s-1},m_s).$$ \par
{\it Then the following sequence (defined via the
multiplication by h) is exact:}
$$0 \rightarrow \left( {I^\prime \over {J^\prime}} \right) _{t-1}
\rightarrow \left( {I \over J} \right) _t
\rightarrow \left( {I + (h)\over J+(h)} \right) _t
\rightarrow 0 . $$
\medskip
\noindent {\it Proof:} Consider the exact sequence (defined
by multiplication by $h$):
$$ 0 \rightarrow {I^\prime \over {J^\prime}}(-1) \rightarrow
{I \over J} \rightarrow {{I \over J}\over Im {I^\prime \over {J^\prime}}}
$$
We have:
$$ Im {I^\prime \over {J^\prime}} = {hI^\prime \over h{J^\prime}} =
{(h)\cap I \over (h)\cap J} =
{(h)\cap I \over (h)\cap J\cap I} = $$
$$= {((h)\cap I ) + J\over J} =  {((h)+J)\cap (J+I )\over J} =
{(((h)+J)\cap I )\over J}. $$
Hence:
$${{I \over J}\over Im {I^\prime \over {J^\prime}}} =
{{I \over J}\over {((h)+J)\cap I \over J}} =
{I \over ((h)+J)\cap I} =
{I+((h)+J)\over (h)+J} = {(h)+I \over (h)+J}.
$$
\qquad  \qquad \qquad \qquad \qquad \qquad \qquad
\qquad  \qquad \qquad \qquad \qquad \qquad \qquad \qed
\medskip

\noindent {\bf Lemma 4.2:} {\it In case 4)  of Remark 3.6,  let }
$$N\prim = (Q,Q_1,...,Q_{s-1};n,(n_1-1)^+,(n_2-1)^+,n_3,...,n_{s-1}).$$
{\it Then } $(I_{N})_{m_s}$ {\it and} $(I_{N\prim })_{m_s-1}$ {\it are
regular.}
\par
\medskip
\noindent {\it Proof:} The conclusion follows by Lemma 3.10 applying a result
by B. Segre (see \ref{6}, \ref{2})
which says that the linear system of curves of degree $d$
in $\PP 2$, with multiplicities $\alpha _1\geq ...\geq \alpha _s$
at points $P_1,...,P_s$ lying on a non singular conic
is regular if and only if
$$d \geq max\{\alpha _1+\alpha _2-1,
\left[{\sum ^{s} _{i=1}\alpha _i\over 2}\right]\}.$$
\qquad  \qquad \qquad \qquad \qquad \qquad \qquad
\qquad  \qquad \qquad \qquad \qquad \qquad \qquad
\qquad \qquad \qquad \qquad \qed
\medskip
\medskip

Now let $J\prim,I\prim $,$h$ be as in Lemma 4.1, and  $N\prim $ as
in the Lemma 4.2. Consider the following diagram:

$$\matrix{&&0&&0&&0&&\cr &&\downarrow&&\downarrow&&
\downarrow&& \cr
0&\rightarrow &\left({I\prim \over J\prim}\right)_{t-1}&
\rightarrow &
\left({I \over J}\right)_{t}
&\rightarrow &\left({I + (h)\over J + (h)}\right)_t&
\rightarrow &0\cr
&&&&&&&& \cr
&&\downarrow&&\downarrow \phi&&\downarrow \psi && \cr
&&&&&&&& \cr
0&\rightarrow &(I_{N\prim })_{m_s-1}&\rightarrow &
(I_N)_{m_s} &\rightarrow &{\bf K}&
\rightarrow &0\cr
&&\downarrow&&\downarrow&&\downarrow&& \cr
&&0&&{\rm coker}\phi&&{\rm coker}\psi&&}$$
\medskip
\noindent
The first exact sequence comes from
Lemma 4.1 and the second is defined via multiplication
by a linear form defining the line $Q_1Q_2$ (${\bf K}$ being its cokernel).
  The first vertical sequence is exact by
induction hypothesis, the map $\phi $ (if
we take $P_s$ to be the origin) comes from the map which
 associates to each $F\in I _t$ the
 tangent cone to $\{ F=0 \}$ (whose equation lies in
 $(I_N)_{m_s}$). We know, by Lemma 3.1, that $\phi $ is injective,
while the map $\psi $ is injective by the Snake's Lemma.
\par
\medskip
{}From Lemma 4.2 it follows (notations as above):
\par
\medskip

\noindent {\bf Corollary 4.3:} {\it In case 4) of Remark 3.6 we have}
$\dim {\bf K} = m_s+1-n_1-n_2${\it .}
 \par
\medskip
{\it Proof:} Since $(I_{N})_{m_s}$ and $(I_{N\prim })_{m_s-1}$  are regular,
 we have \par
\noindent
$\dim {\bf K} = \dim (I_{N})_{m_s} - (I_{N\prim })_{m_s-1} = $ \par

$={m_s+2\choose 2} - \sum ^{s-1} _{i=1}{n_i+1\choose 2} - {n+1\choose 2}
 - {m_s+1\choose 2} + \sum ^{s-1} _{i=1}{n_i+1\choose 2} - n_1 - n_2
 + {n+1\choose 2} =$  \par
\centerline{$m_s + 1 -n_1 - n_2 .$}
\par
\qquad \qquad \qquad \qquad
\qquad \qquad \qquad \qquad
\qquad \qquad \qquad \qquad
\qquad \qquad  \qed
\medskip
\medskip
Thus if we prove that, in case $4),$
$\dim \left( {I + (h)\over J + (h)}\right) _t \geq m_s + 1-n_1-n_2 ,$
 then the map $\psi$
will be surjective, and we will be done ($\phi$ will be surjective too,
hence an isomorphism).
\par
First we deal with a particular case of 4), namely $n_1=0$.
 We will prove \medskip

\noindent {\bf Lemma 4.4:} {\it Let } $s \geq 5$,
$ n \leq m_s$, $m_s>0$ {\it and} $n_1 = 0$ (i.e. $t \leq m_1+m_s).$
{\it Then }
$$\dim \left( {I + (h)\over J + (h)}\right) _t \geq m_s + 1 .$$
\medskip
\noindent {\it Proof:} We have to find $m_s+1$ linearly independent
forms in $I$ such that their classes remain independent
modulo $(J + (h))\cap I$ (see proof of Lemma 4.1). \par
\medskip
Let $G_1$,$G_2$ be the quadric forms defining the cones on $\Gamma $
with vertices in $P_1$,$P_2$ respectively.  \par
Let $m_1=\ldots = m_s$. Then it is easy to check that the forms
$G_1 ^xG _2 ^yH^{t-2m_s}$
give what we want, when $H$ is a plane not containing any of the
$P_i$'s and $x+y = m_s$.
\par
\medskip
If $m_s = 1$, $m_1 > m_s$, we consider a form
$S \in ({\goth p}^{m_1-2} _1 \cap {\goth p}^{m_2-1} _2\cap \ldots
\cap {\goth p}^{m_{s-1}-1} _{s-1})_{t-2}$ which is not zero at $P_s$.
In order to prove that such $S$ exists, the following Lemma
(from [{\bf 3}], Lemma 4) will be useful: \par
\medskip

\noindent {\bf Lemma 4.5:} {\it Let }
$P_1,...,P_\lambda,P \in \PP r$ {\it be distinct points
in general (linear) position and let}
$Y={\goth p}^{m_1} _1 \cap \ldots \cap {\goth p}^{m_{\lambda}} _\lambda,$
{\it with} $m_1\geq ... \geq m_\lambda \geq 0$. {\it If} $t\in {\Bbb N}$
{\it is such that} $rt \geq \sum ^\lambda _{i=1} m_i$ {\it and}
$t\geq m_1${\it , then we can find a form} $F\in Y_t$
{\it avoiding P}. \par
\medskip
So, by this lemma, to check that $S$ exists,
it suffices to show that \par
 \qquad $3(t-2)\geq \sum ^s _{i=1} m_i -(s+1)$,
$t-2\geq m_1-2$,\  and \  $t-2\geq m_2-1$.
\par
The first inequality is verified, since $n\leq m_s =1$ and
$s \geq 5$ imply
that ${\cal P}(2)$ does not hold for $I_t$,
i.e. $3t-5 \geq \sum ^s _{i=1} m_i -s$.
The other two because $n_1=0$, hence
$t\geq m_1+m_s=m_1+1$.
\par
In the same way one can find a form
$S\prim \in ({\goth p}^{m_1-1} _1 \cap {\goth p}^{(m_2-2)^+} _2
\cap \ldots
\cap {\goth p}^{m_{s-1}-1} _{s-1})_{t-2}$ which is not zero at $P_s$.
We only notice that for $m_2=1$, we have to check that $3(t-2) \geq
\sum ^s _{i=1} m_i -s = m_1-1$, and this follows from
$t \geq m_1+1$.
\par
Then the forms  $G_1 S$, and $G_2 S\prim $  will give
what we want.
\medskip
Now let $m_s > 1$ and $m_1 > m_s$.  We can work by induction
on $\sum ^s _{i=1} m_i$, using the previous cases as initial steps).
\par
If $m_1>m_2$, consider the homogeneous ideals  $I^*$,$J^*$
associated to the schemes
$$Z^*=(P_1,...,P_s;m_1-2,m_2-1,m_3-1,...,m_s-1) =
(P_1,...,P_s;m_1^*,m_2^*,...,m_s^*),$$
$$W^*=(P_1,...,P_s;m_1-2,m_2-1,...,m_{s-1}-1,m_s)=
(P_1,...,P_s;m_1^*,m_2^*,...,m_s^*+1).$$

By induction we get that there exist
forms $F_1,...,F_{m_s} \in I^* + (h)$ such that their
classes in ${I^* + (h)\over J^* + (h)}$ are linearly
independent. Let us check that
$t-2$,$I^*$,$J^*$ verify the hypotheses of the Lemma;
since  ${\cal P}(m_s+1)$ doesn't hold
 for $I_t$, then \par
\centerline{$3t-5m_s \geq \sum ^s _{i=1} m_i -sm_s$.}
\par
 It follows \
$3(t-2)+5(1-m_s) \geq \sum ^s _{i=1} m_i -sm_s-1 =
\sum ^s _{i=1} m_i^* +s-sm_s $, so (with obvious notation)\
$n^* \leq m_s-1 = m_s^*$.\ Moreover $n_1^*=0$,
in fact $t-2 \geq m_1+m_s-2 \geq m_1^*+m_s^*$, and
$m_1^* \geq m_2^* \geq ... \geq m_s^*>0$, since
$m_1  > m_2  \geq ... \geq m_s >1$.
\par
Hence we can consider the forms $G_1F_1,...,G_1F_{m_s}$
which are independent modulo $(J+(h))\cap I $, as required. \par
We have to find another one, so consider the number:
$A = \sum ^s _{i=1} m_i - m_1 - (s-1)m_s $.
When $A \geq m_s$, consider the quadric forms
$G _2$,$\ldots $,$G _{s-1}$ defining the cones on $\Gamma $
with vertices in
$P_2$,...,$P_{s-1}$ respectively, and let
$F =  G _2^{m_2-m_s}G _3^{m_3-m_s}\ldots $, where we go
on with the products until we get $m_s$ factors (i.e.
$\deg F = 2m_s$). Let $F \in
{\goth p}^{m_s} _1 \cap {\goth p}^{m\prim _2}_2\cap \ldots
\cap {\goth p}^{m\prim _{s-1}} _{s-1}\cap {\goth p}^{m_s} _s$,
then we can choose another form $S \in
({\goth p}^{m_1-m_s} _1 \cap {\goth p}^{m _2-m\prim_2}_2 \cap \ldots
\cap {\goth p}^{m _{s-1}-m\prim_{s-1}} _{s-1})_{t-2m_s}$,
not passing through
$P_s$ (this is possible by Lemma 4.5
since $t-2m_s\geq m_1-m_s$, and $3(t-2m_s)\geq \sum ^s _{i=1} m_i -sm_s-m_s$).
\par
The forms $FS,G_1F_1,...,G_1F_{m_s}$ are what we want.
\par
When $A<m_s$ we consider
$F = G _2^{m_2-A}G _3^{m_3-m_s}\ldots
G _{s-1}^{m_{s-1}-m_s}$  and
$S \in {\goth p}^{m_1-m_s} _1$ not passing through $P_s $
(again, it exists since $t-2m_s\geq m_1-m_s$ and $3(t-2m_s)\geq
m_1-m_s$). \par
$FS,G_1F_1,...,G_1F_{m_s}$ are the forms that
we want.
\medskip
Finally, when $m_s>1$, $m_1=m_2=\ldots =m_r>m_{r+1}$, $r<s$,
we proceed as before, but starting with
ideals $I^*$,$J^*$ associated to the schemes:
$$Z^*=(P_1,...,P_s;m_1-1,...,m_{r-1}-1,m_r-2,m_{r+1}-1,...,m_s-1)
=(P_1,...,P_s;m_1^*,...,m_s^*),$$
$$W^*=(P_1,...,P_s;m_1-1,...,m_{r-1}-1,m_r-2,m_{r+1}-1,...,m_{s-1}-1,m_s),$$
so that $m_1^* \geq m_2^* \geq ... \geq m_s^*>0$,\   $\{h^*=0\}$
is the plane $P_1P_2P_s$
 again, and we use $G _r$ instead of
$G _1$.   \par
\qquad  \qquad \qquad \qquad \qquad \qquad \qquad
\qquad  \qquad \qquad \qquad \qquad \qquad \qquad
\qed
\medskip
\medskip
With the following proposition the proof of Theorem 2.2 will
be complete (by Corollary 4.3).
\par
\medskip
\medskip

\noindent {\bf Proposition 4.6:} {\it In case 4) of Remark 3.6
we have:}
$$\dim \left({I + (h)\over J + (h)}\right)_t \geq m_s + 1 -n_1 -n_2 .$$
{\it Hence the map } $\psi $ is surjective.
\medskip
{\it Proof:} If $n_1=0$ the conclusion follows by Lemma 4.4.
\par
Assume $n_1>0.$
Let us consider the form
$S = G^{n_1-n_2}_1.G^{n_2}_s$, where, as usual,
$G_i$ defines the cone on $\Gamma$ with vertex in $P_i$.
\par
The form $S$ has degree $2n_1$, multiplicity $2n_1-n_2$ at
$P_1,$ $n_1$ at $P_2,...,P_{s-1}$ and $n_1+n_2$ at
$P_s.$ \par
Now let $$I^*= {\goth p}^{m_1-2n_1+n_2} _1
\cap {\goth p}^{m_2-n_1}_2 \cap \ldots
\cap {\goth p}^{m_{s-1}-n_1} _{s-1}\cap
{\goth p}^{m_s-n_1-n_2} _s =
{\goth p}^{m_1^*} _1 \cap \ldots
\cap {\goth p}^{m_s^*} _s $$ and
$$J^* = {\goth p}^{m_1-2n_1+n_2} _1
\cap {\goth p}^{m_2-n_1}_2 \cap \ldots
\cap {\goth p}^{m_{s-1}-n_1} _{s-1}\cap
{\goth p}^{m_s-n_1-n_2+1} _s.$$
\par
It is easy to check that $m_1^* \geq m_2^* \geq \ldots \geq m_s^* \geq 0.$ \par
Since $t \geq 2n_1$ and any base of
$ \left({I^* + (h)\over J^* + (h)}\right)_{t-2n_1}$, when
multiplied by $S$, gives independent elements
in $  \left({I + (h)\over J + (h)}\right)_{t-2n_1},$ we have
$$\dim \left({I + (h)\over J + (h)}\right)_t \geq
\dim \left({I^* + (h)\over J^* + (h)}\right)_{t-2n_1}.$$
Hence we will be finished if we show that
\par
\medskip
\noindent $(5)  \qquad \qquad \qquad \qquad \qquad
\dim \left({I^* + (h)\over J^* + (h)}\right)_{t-2n_1}
\geq m_s-n_1-n_2+1.$
\medskip
Let us consider first the case when $m_s=n_1+n_2$. In this situation we just
have to find one form $F$ in $I^*_{t-2n_1}$
which is not zero at $P_s.$
\par
\medskip
By Lemma 4.5 it suffices to check that :
\par
\medskip
\noindent a) $3(t-2n_1)\geq \sum ^s _{i=1} m_i - (s+1)n_1$; \par
\noindent b) $t-2n_1 \geq m_1-2n_1+n_2$.\par
\medskip
Since $n \leq m_s-n_1 = n_2 \leq n_1,$ then
 ${\cal P}(n_1+1)$ does not hold for $I_t$, i.e.
$3t-5n_1 \geq \sum ^s _{i=1} m_i -s n_1$, hence a) holds.\par
 From $m_s=n_1+n_2$ we get $t = m_1+m_s-n_1 = m_1+n_2$, which is b).
\par
Now let $m_s>n_1+n_2.$  If we can apply Lemma 4.4 to
$I^* , J^*, h,t-2n_1,$ we get exactly
the inequality (5) and we will be done. Thus let us check
that the hypotheses of Lemma 4.4 are satisfied.
\par
We have $m^*=m_s-n_1-n_2>0$, so the number of non-zero exponents in
$I^* $ is exactly $s\geq 5$.
\par
We also have $t-2n_1 = (m_1+m_s-n_1)-2n_1=m_1^*+m_s^*.$ \par
Finally since $ n \leq m_s,$ we have that ${\cal P}(n+1)$ does not
hold for $I_t,$ i.e.
$3t-5n \geq \sum ^s _{i=1} m_i -sn$.
  This implies that
$3(t-2n_1)-5(n-n_1) \geq \sum ^s _{i=1} m_i^* -s(n-n_1)$.  Hence,
for $n \geq n_1,$ we get (with obvious notation)
$n^* \leq n-n_1 \leq m_s-n_1-n_1 \leq m_s-n_1-n_2 = m_s^*$.
For $n < n_1,$ since $s \geq 5$, we have $n^* = 0 \leq m_s^*.$ \par
Thus in both cases we get $n^* \leq m_s^*$ , so all the hypotheses of Lemma 4.4
are satisfied,
and the proof is complete.

\qquad  \qquad \qquad \qquad \qquad \qquad \qquad
\qquad  \qquad \qquad \qquad \qquad \qquad \qquad
\qed
\medskip
\medskip
\centerline{{\bf REFERENCES}}
\medskip
\medskip
\noindent [{\bf 1}]: R.Achilles, P.Schenzel, W.Vogel: {\it
Bemerkungen uber Normale Flachheit und\par
\noindent Normale
Torsionfreiheit und Anwendungen}. Per. Mat. Hung. {\bf 12},
 (1981), 49-75. \par
\medskip
\noindent [{\bf 2}]: M.V.Catalisano:
{\it "Fat" points on a conic}. Comm. Algebra {\bf 19}(8), (1991), 2513-2168.
\par
\medskip
\noindent [{\bf 3}]: M.V.Catalisano, Trung, G.Valla:
{\it A sharp bound for the regularity index of fat points in general
 position.} Proc. Amer. Math. Soc. {\bf 118} (1993), 717-724.
\par
\medskip
\noindent [{\bf 4}]: J.Harris: {\it Algebraic Geometry}.
 Springer-Verlag, Grad. Texts in Math. {\bf 133} (1992), Berlin,
 New York.\par
\medskip
\par
\noindent [{\bf 5}]: L.Robbiano:
{\it An algebraic property} of $
\PP 1 \times \PP N$. Comm.in Alg. {\bf 7}, (1979),
641-655.
 \par
\medskip
\noindent \ref{6}: B.Segre: {\it Alcune questioni su insiemi
finiti di punti in Geometria Algebrica.} Atti Convegno
Internazionale di Geometria Algebrica. Torino (1961), 15-33.
\par
\medskip
\noindent \ref{7}: N.V.Trung, G.Valla: {\it Upper bounds for the
regularity index of fat points with uniform position property}.
J. Algebra, to appear.

\medskip
\medskip
\medskip
\medskip
\medskip
\medskip
\medskip
\medskip
Address of the authors:
M.V.Catalisano: Dip. di Matematica, Univ. di Genova, Via L.B.Alberti 4, Genova,
Italy. {\it e-mail}: catalisa@dima.unige.it
\medskip
A.Gimigliano: Dip. di Matematica, Univ. di Bologna, P.zza di Porta S.Donato 5,
Bologna, Italy. {\it e-mail}: gimiglia@dm.unibo.it

 \end